\documentclass[12pt,thmsa]{article}%
\usepackage{amsfonts}
\usepackage{sw20bams}
\usepackage{amsmath}
\usepackage{amssymb}
\usepackage{graphicx}%
\providecommand{\U}[1]{\protect\rule{.1in}{.1in}}
\begin{document}

\author{Steven R. Finch and John A. Shonder}
\title{Lost at Sea}
\date{March 12, 2016}
\maketitle

\begin{abstract}
What is the path of minimum expected length for escaping a planar convex
region $\Omega$? We rigorously obtain best $2$-segment and $3$-segment
solutions when $\Omega$ is an infinite strip, and numerically examine
$2$-segment solutions when $\Omega$ is a disk.

\end{abstract}

\footnotetext{Copyright \copyright \ 2004, 2016 by Steven R. Finch. All rights
reserved.}A swimmer is lost in a dense fog at sea. She knows that the sea is a
planar infinite strip
\[
\left\{  (x,y)\in\mathbb{R}^{2}:0\leq x\leq1\right\}
\]
of unit width. Assume that the $x$-coordinate of her initial position is
uniformly distributed on the interval $[0,1]$; the distribution of the
$y$-coordinate can be arbitrary. Assume as well that her initial orientation
(the angle between her initial velocity vector and the $x$-axis) is uniformly
distributed on $[-\pi,\pi]$ and that her speed is constant. What escape
trajectory should the swimmer follow that minimizes her expected time to reach
either shore?

The min-max analog of this problem was solved long ago (see \cite{FW} for a
survey). Zalgaller \cite{Za1, Za2} proposed a heuristic solution of the above
min-mean problem. Shonder \cite{Sh} independently determined the precise
$2$-segment solution of least expected escape time. (We agree that a
$k$-segment path is a continuous, piecewise linear curve consisting of $\leq
k$ pieces.)\ \ The details underlying Shonder's computation appear here for
the first time. We find a $3$-segment path that improves slightly upon the $2
$-segment path and also demonstrate that Zalgaller's solution is far from optimal.

Next, we examine a different sea:\ the disk
\[
\left\{  (x,y)\in\mathbb{R}^{2}:x^{2}+y^{2}\leq1\right\}
\]
of unit radius. Assume that the starting point is uniformly distributed on the
disk and that we seek (as before) to minimize the expected time to reach the
boundary. Intuition suggests that the swimmer should follow a $1$-segment
path. If this is true, then the mean length of an arbitrary escape trajectory
is $\geq8/(3\pi)=0.8488263631...$ \cite{Co, Cu}. A calculus-of-variations
proof of this general inequality is not known. We can confirm this only
numerically for $2$-segment paths.

In the final section, there appears a first attempt at proving the above
inequality under special circumstances. The technique is due to Gevirtz
\cite{Ge} and yields results when the escape path curvature is sufficiently small.

\subsection{Infinite Strip:\ 2-Segment Scenario}

Without loss of generality, let the initial position be $(x,0)$ and the
initial orientation be $\theta$, where $0\leq x\leq1$ and $0\leq\theta\leq\pi
$. Fix a distance $r\geq0$ and an angle $0\leq\alpha\leq\pi$. As shown in
Figure 1, there are five distinct cases to consider.\newline%
%TCIMACRO{\FRAME{ftbpFU}{6.0658in}{2.4656in}{0pt}{\Qcb{Five cases for the
%infinite strip.}}{}{las.eps}{\special{ language "Scientific Word";
%type "GRAPHIC";  maintain-aspect-ratio TRUE;  display "USEDEF";
%valid_file "F";  width 6.0658in;  height 2.4656in;  depth 0pt;
%original-width 7.9321in;  original-height 3.2085in;  cropleft "0";
%filename 'las.eps';file-properties "XNPEU";}} }%
%BeginExpansion
\begin{figure}[ptb]%
\centering
\includegraphics[
height=2.4656in,
width=6.0658in
]%
{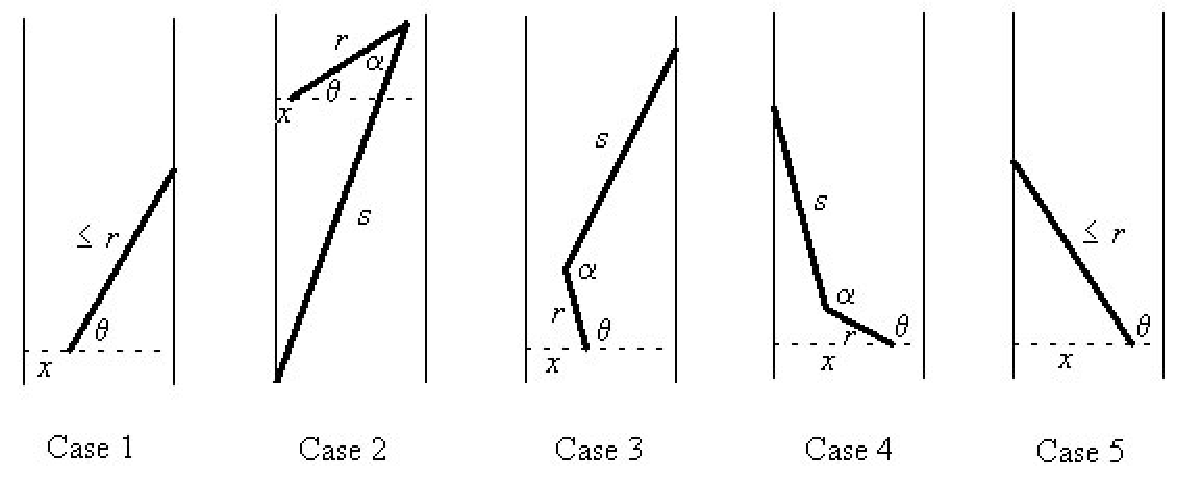}%
\caption{Five cases for the infinite strip.}%
\end{figure}
%EndExpansion

\textbf{Case 1}:\ the swimmer reaches the right-hand shore before traveling a
distance $r$, which happens when $x+r\geq1$ and $\theta\leq\arccos((1-x)/r)$.
The path length is $(1-x)/\cos(\theta)$.

\textbf{Case 2}: the swimmer travels a distance $r$, does not reach the
right-hand shore, pivots an angle $\alpha$, and then travels until the
left-hand shore is reached. This happens when
\[
\left\{  \left[  x+r\geq1\text{ and }\arccos((1-x)/r)<\theta\right]  \text{ or
}\left[  x+r<1\right]  \right\}  \,\text{ and }\left\{  \theta<\tfrac
\pi2-\alpha\right\}  .
\]
The path length is $r+\left(  x+r\cos(\theta)\right)  /\cos(\theta+\alpha) $.

\textbf{Case 3}:\ the swimmer travels a distance $r$, does not reach a shore,
pivots an angle $\alpha$, and then travels until the right-hand shore is
reached. This happens when
\begin{align*}
&  \left\{  \left[  x+r\geq1\text{ and }x-r\leq0\text{ and }\arccos
((1-x)/r)<\theta<\pi-\arccos(x/r)\right]  \right. \\
&  \text{or }\left[  x+r<1\text{ and }x-r\leq0\text{ and }\theta<\pi
-\arccos(x/r)\right]  \text{ }\\
&  \text{or }\left[  x+r\geq1\text{ and }x-r>0\text{ and }\arccos
((1-x)/r)<\theta\right]  \text{ }\\
&  \left.  \text{or }\left[  x+r<1\text{ and }x-r>0\right]  \right\}  \text{
and }\left\{  \tfrac\pi2-\alpha\leq\theta\leq\tfrac{3\pi}2-\alpha\right\}  .
\end{align*}
The path length is $r+\left(  x-1+r\cos(\theta)\right)  /\cos(\theta+\alpha)$.

\textbf{Case 4}: the swimmer travels a distance $r$, does not reach the
left-hand shore, pivots an angle $\alpha$, and then travels until the
left-hand shore is reached. This happens when
\[
\left\{  \left[  x-r\leq0\text{ and }\theta<\pi-\arccos(x/r)\right]  \text{ or
}\left[  x-r>0\right]  \right\}  \,\text{ and }\left\{  \tfrac{3\pi}%
2-\alpha<\theta\right\}  .
\]
The path length is $r+\left(  x+r\cos(\theta)\right)  /\cos(\theta+\alpha) $.

\textbf{Case 5}:\ the swimmer reaches the left-hand shore before traveling a
distance $r$, which happens when $x-r\leq0$ and $\pi-\arccos(x/r)\leq\theta$.
The path length is $-x/\cos(\theta)$.\newline

We've restricted attention to $0\leq\theta\leq\pi$ for simplicity's sake. The
set of escape paths corresponding to $-\pi\leq\theta\leq0$ is obtained from
the set of escape paths corresponding to $0\leq\theta\leq\pi$ via reflection
across the horizontal axis, followed by reflection across the vertical axis.
This composite transformation can be written as $(x,\theta)\mapsto
(1-x,\theta-\pi)$.

If it is assumed that $r>1$ and $0<\frac{\pi}{2}-\alpha<\arccos(1/r)$, then
only Cases 1, 3 and 5 enter into the calculations of the expected path
length:
\[
\frac{1}{\pi}%
%TCIMACRO{\dint \limits_{0}^{1}}%
%BeginExpansion
{\displaystyle\int\limits_{0}^{1}}
%EndExpansion
\left[  \int\limits_{0}^{\arccos\left(  \frac{1-x}{r}\right)  }\frac{1-x}%
{\cos(\theta)}d\theta+\int\limits_{\arccos\left(  \frac{1-x}{r}\right)  }%
^{\pi-\arccos\left(  \frac{x}{r}\right)  }\left(  r+\frac{x-1+r\cos(\theta
)}{\cos(\theta+\alpha)}\right)  d\theta-\int\limits_{\pi-\arccos\left(
\frac{x}{r}\right)  }^{\pi}\frac{x}{\cos(\theta)}d\theta\right]  dx.
\]
All of these integrals can be evaluated in closed form. Hence the expected
path length, multiplied by $\pi$, becomes
\begin{align*}
&  \ \ \ \ r\left[  \left(  1+\cos\alpha\right)  \left(  \frac{\pi}{2}%
+\sqrt{r^{2}-1}-r\right)  -\left(  2+\cos\alpha\right)  \arccos\left(
\frac{1}{r}\right)  +\frac{\pi}{2}\right]  +\\
&  \ \ \ \ \ln\left(  \sqrt{r^{2}-1}+r\right)  +\frac{1}{2}r^{2}\ln
(1-\cos\alpha)\sin^{2}\alpha+r\sin\alpha\ln\left(  \frac{\sqrt{r^{2}-1}%
\sin\alpha-\cos\alpha}{r\sin\alpha}\right)  +\\
&  \ \ \ \ \frac{1}{4}\left(  r^{2}\sin^{2}\alpha+1\right)  \ln\left(
\sqrt{r^{2}-1}\cos\alpha+\sin\alpha+r\right)  +\\
&  \ \ \ \ \frac{1}{4}\left(  r^{2}\sin^{2}\alpha-1\right)  \ln\left(
\sqrt{r^{2}-1}\cos\alpha-\sin\alpha+r\right)  -\\
&  \ \ \ \ \frac{1}{4}\left(  r^{2}\sin^{2}\alpha+1\right)  \ln\left[
-\sqrt{r^{2}-1}\cos^{2}\alpha-\left(  \sin\alpha+\sqrt{r^{2}-1}-r\right)
\cos\alpha-\sin\alpha+r\right]  -\\
&  \ \ \ \ \frac{1}{4}\left(  r^{2}\sin^{2}\alpha-1\right)  \ln\left[
-\sqrt{r^{2}-1}\cos^{2}\alpha+\left(  \sin\alpha-\sqrt{r^{2}-1}+r\right)
\cos\alpha+\sin\alpha+r\right]  .
\end{align*}
Minimizing this expression gives $r=1.0432668686...$ and $\alpha
=1.3734935859...\approx78.7^{\circ}$ as the optimal parameter values for the
$2$-segment scenario. Therefore the least expected path length is
$0.8869669056...$. See Figure 2 for several sample realizations (with fixed
$r$ and $\alpha$).%
%TCIMACRO{\FRAME{ftbpFU}{3.2534in}{5.2062in}{0pt}{\Qcb{Four sample realizations
%with $r=1.043$, $\alpha=78.7^{\circ}$, $x=0.1,0.3,0.5,0.9$ and $\theta
%=80^{\circ},160^{\circ},50^{\circ},140^{\circ}$ (involving three Cases).}}%
%{}{pth2.eps}{\special{ language "Scientific Word";  type "GRAPHIC";
%width 3.2534in;  height 5.2062in;  depth 0pt;  original-width 3.2085in;
%original-height 5.1508in;  cropleft "0";  croptop "1";  cropright "1";
%cropbottom "0";  filename 'pth2.eps';file-properties "XNPEU";}}}%
%BeginExpansion
\begin{figure}[ptb]%
\centering
\includegraphics[
height=5.2062in,
width=3.2534in
]%
{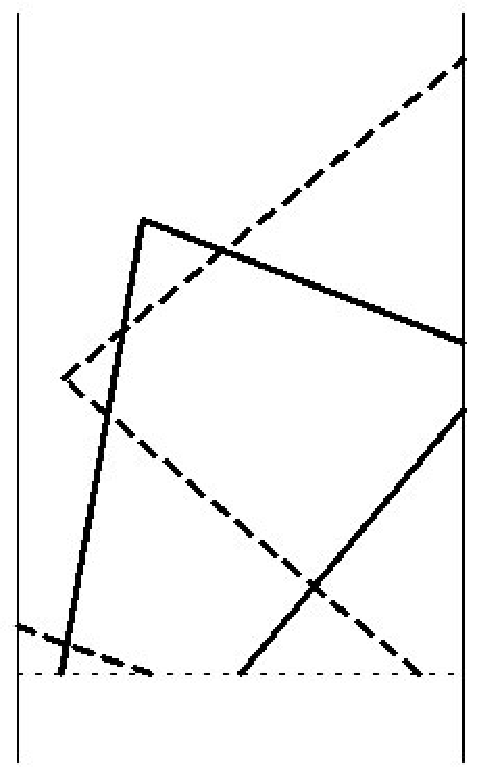}%
\caption{Four sample realizations with $r=1.043$, $\alpha=78.7^{\circ}$,
$x=0.1,0.3,0.5,0.9$ and $\theta=80^{\circ},160^{\circ},50^{\circ},140^{\circ}$
(involving three Cases).}%
\end{figure}
%EndExpansion

Here is a new problem. For each $(r,\alpha)$, there exists a probability
distribution of escape times obtained by sampling $(x,\theta)$ uniformly from
$[0,1]\times[0,\pi]$. Compute the median of each distribution and then
determine $(r,\alpha)$ corresponding to the least median. Monte Carlo
simulation suggests that a solution is $r=\infty$ ($\alpha$ can be arbitrary),
consistent with the idea that medians are less sensitive to outliers than
means. Define a random variable $q$ to be $(1-x)/\cos(\theta) $ if $\theta
<\pi/2$ and $-x/\cos(\theta)$ if $\theta>\pi/2$. The path length is the median
of $q$, which is estimated to be $0.78...$.

\subsection{Infinite Strip:\ 3-Segment Scenario}

We assume that $r>1$ and $0<\frac{\pi}{2}-\alpha<\arccos(1/r)$ as at the end
of the preceding section. Cases 1 and 5 thus remain unchanged and only Case 3
needs to be refined. Fix a second distance $s\geq0$ and a second angle
$0\leq\beta\leq\pi$. See Figure 3.\newline%
%TCIMACRO{\FRAME{ftbpFU}{5.5988in}{4.1399in}{0pt}{\Qcb{Two subcases of Case
%3.}}{}{las2.eps}{\special{ language "Scientific Word";  type "GRAPHIC";
%width 5.5988in;  height 4.1399in;  depth 0pt;  original-width 5.5417in;
%original-height 4.0906in;  cropleft "0";  croptop "1";  cropright "1";
%cropbottom "0";  filename 'las2.eps';file-properties "XNPEU";}} }%
%BeginExpansion
\begin{figure}[ptb]%
\centering
\includegraphics[
height=4.1399in,
width=5.5988in
]%
{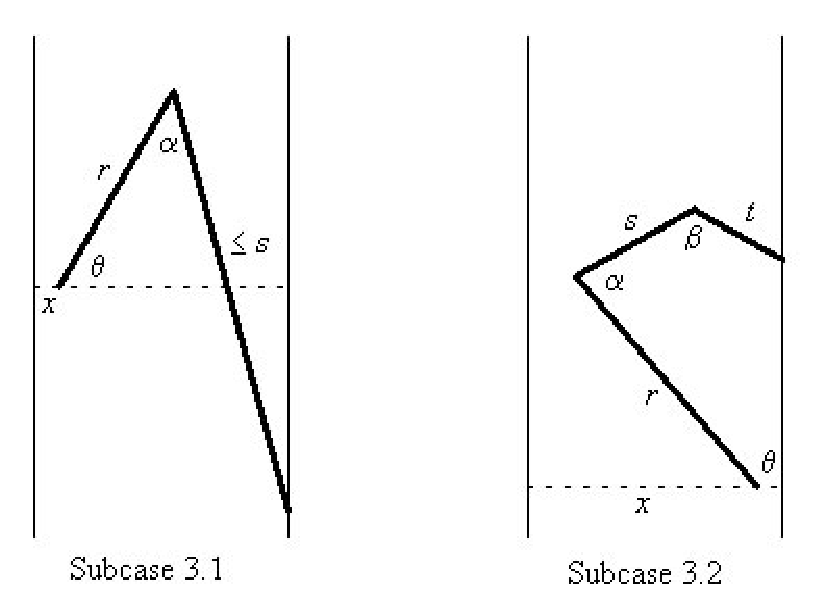}%
\caption{Two subcases of Case 3.}%
\end{figure}
%EndExpansion

\textbf{Subcase 3.1}:\ the swimmer travels a distance $r$, does not reach a
shore, pivots an angle $\alpha$, and then reaches the right-hand shore before
traveling a distance $s$. This happens when
\[
\arccos((1-x)/r)<\theta\leq\theta_{0}(x,r,s,\alpha)
\]
where $\theta_{0}$ is a solution of the equation
\[
1-x-r\cos(\theta)+s\cos(\theta+\alpha)=0;
\]
an exact formula for $\theta_{0}$ is
\[
\theta_{0}=\arcsin\left(  \frac{s\sin(\alpha)}{\sqrt{r^{2}+s^{2}%
-2\,r\,s\cos(\alpha)}}\right)  +\arccos\left(  \frac{1-x}{\sqrt{r^{2}%
+s^{2}-2\,r\,s\cos(\alpha)}}\right)  .
\]
The path length is $r+\left(  x-1+r\cos(\theta)\right)  /\cos(\theta+\alpha)$.

\textbf{Subcase 3.2}:\ the swimmer travels a distance $r$, does not reach a
shore, pivots an angle $\alpha$, travels a distance $s$, does not reach a
shore, pivots an angle $\beta$, and then travels until the right-hand shore is
reached. This happens when
\[
\theta_{0}(x,r,s,\alpha)<\theta<\pi-\arccos(x/r)
\]
and $\beta$ satisfies other constraints (for example, if $\beta$ is too small,
then the $3$-segment might cross itself and reach the left-hand shore rather
than the right-hand shore). We will not attempt to prescribe these additional
inequalities. The path length is $r+s+\left(  1-x+s\cos(\theta+\alpha
)-r\cos(\theta)\right)  /\cos(\theta+\alpha+\beta).$\newline

Under the proper circumstances, therefore, the expected path length is
\begin{align*}
&  \ \ \ \frac1\pi%
%TCIMACRO{\dint \limits_{0}^{1}}%
%BeginExpansion
{\displaystyle\int\limits_{0}^{1}}
%EndExpansion
\left[  \int\limits_{0}^{\arccos\left(  \frac{1-x}r\right)  }\frac{1-x}%
{\cos(\theta)}d\theta+\int\limits_{\arccos\left(  \frac{1-x}r\right)
}^{\theta_{0}}\left(  r+\frac{x-1+r\cos(\theta)}{\cos(\theta+\alpha)}\right)
d\theta+\right. \\
&  \ \ \ \left.  \int\limits_{\theta_{0}}^{\pi-\arccos\left(  \frac xr\right)
}\left(  r+s+\frac{1-x+s\cos(\theta+\alpha)-r\cos(\theta)}{\cos(\theta
+\alpha+\beta)}\right)  d\theta-\int\limits_{\pi-\arccos\left(  \frac
xr\right)  }^{\pi}\frac x{\cos(\theta)}d\theta\right]  dx.
\end{align*}
A closed-form expression is unlikely here because of the presence of
$\theta_{0}$ in the middle two integrals:\ only integration with respect to
$\theta$ seems feasible. Let
\[%
\begin{array}
[c]{ccc}%
\kappa=\sqrt{r^{2}+s^{2}-2\,r\,s\cos(\alpha)}, &  & \rho=\arcsin\left(
s\sin(\alpha)/\kappa\right)  ,
\end{array}
\]
\[%
\begin{array}
[c]{ccc}%
u_{1}=r\sin(\alpha), &  & u_{2}=s\sin(\beta)-r\sin(\alpha+\beta),
\end{array}
\]
\[%
\begin{array}
[c]{ccc}%
v_{1}=r(1+\cos(\alpha)), &  & v_{2}=r+s+s\cos(\beta)-r\cos(\alpha+\beta),
\end{array}
\]
\[%
\begin{array}
[c]{ccc}%
\xi_{1}=\alpha+\arccos\left(  (1-x)/r\right)  , &  & \xi_{2}=\alpha+\beta
+\rho+\arccos\left(  (1-x)/\kappa\right)  ,
\end{array}
\]
\[%
\begin{array}
[c]{ccc}%
\eta_{1}=\alpha+\rho+\arccos\left(  (1-x)/\kappa\right)  , &  & \eta
_{2}=\alpha+\beta+\pi-\arccos(x/r)
\end{array}
\]
and, for $j=1,2$, define
\[
p_{j}=u_{j}\ln\left|  \frac{\cos(\xi_{j})}{\cos(\eta_{j})}\right|
+\frac{(-1)^{j}(1-x)}2\ln\left(  \frac{\sin(\xi_{j})\left(  \sin(\eta
_{j})+1\right)  -\sin(\eta_{j})-1}{\sin(\xi_{j})\left(  \sin(\eta
_{j})-1\right)  +\sin(\eta_{j})-1}\right)  +v_{j}(\eta_{j}-\xi_{j}).
\]
Then the expected path length, multiplied by $\pi$, becomes
\[
\ln\left(  \sqrt{r^{2}-1}+r\right)  +r\left(  r-\sqrt{r^{2}-1}\right)  +%
%TCIMACRO{\dint \limits_{0}^{1}}%
%BeginExpansion
{\displaystyle\int\limits_{0}^{1}}
%EndExpansion
\left(  p_{1}+p_{2}\right)  dx.
\]
This can be minimized by use of Leibniz's rule for differentiation under the
integral sign. We obtain $r=1.0255050653...$, $\alpha=1.4909825316...\approx
85.4^{\circ}$, $s=0.5306340577...$ and $\beta=2.7495709960...\approx
157.5^{\circ}$ as the optimal parameter values for the $3$-segment scenario,
with least expected path length $0.8835534788...$.

A third subcase, for which $\pi\leq\beta\leq2\pi$, awaits examination. Either
shore could be the final resting place.

\subsection{Infinite Strip:\ Zalgaller's Proposed Solution}

Consider the path in Figure 4, where the points $A=(0,0)$, $B=(0.814,0)$,
$C=(0.8460,0.0005)$, $D=(1.3017,0.0151)$, $E=(0.814,1)$ are joined by line
segments, with the exception of points $B$, $C$ which are joined by a tiny
circular arc of radius $1$, center $E$ and angle $0.032.$ Zalgaller \cite{Za1,
Za2} claimed that, if the swimmer follows this path, then her expected time to
reach either shore is approximately minimal. We examined his claim by
computing best $2$-segment and $3$-segment fits to Zalgaller's path, and then
calculating the expected escape time via our formulas. For example, in the
$2$-segment fit, we obtained $r=1.3017$, $\alpha=64.3^{\circ}$ and hence the
expected escape time is $0.9188$. This is somewhat consistent with Zalgaller's
estimate $0.9523$. More importantly, however, Shonder's optimal $2$-segment
path \cite{Sh} has expected escape time $0.8870$ which is considerably less
than $0.9523$. There must be an error somewhere in the details of \cite{Za1}.
We aim someday to better understand Zalgaller's elaborate construction, in the
hope that his procedure (once corrected) will lead to escape trajectories that
outperform even our optimal $3$-segment path.%
%TCIMACRO{\FRAME{ftbpFU}{3.7421in}{3.672in}{0pt}{\Qcb{Zalgaller's escape
%path.}}{}{zagl.eps}{\special{ language "Scientific Word";  type "GRAPHIC";
%width 3.7421in;  height 3.672in;  depth 0pt;  original-width 3.6945in;
%original-height 3.6244in;  cropleft "0";  croptop "1";  cropright "1";
%cropbottom "0";  filename 'zagl.eps';file-properties "XNPEU";}} }%
%BeginExpansion
\begin{figure}[ptb]%
\centering
\includegraphics[
height=3.672in,
width=3.7421in
]%
{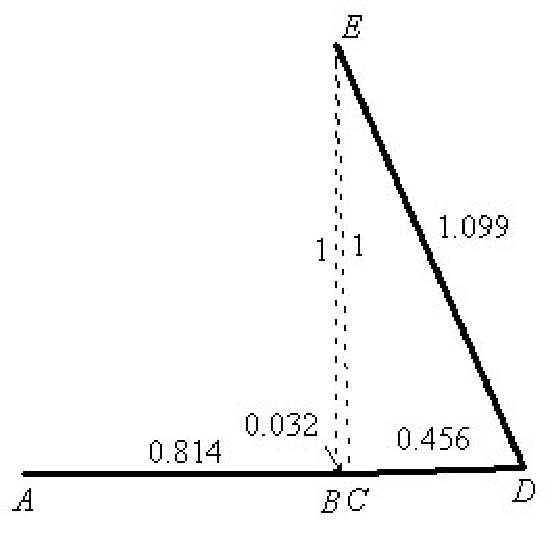}%
\caption{Zalgaller's escape path.}%
\end{figure}
%EndExpansion

\subsection{Disk:\ 2-Segment Scenario}

Without loss of generality, let the starting point be $(x,0)$ and the initial
orientation be $\theta$, where $0\leq x\leq1$ and $-\pi\leq\theta\leq\pi$. Fix
a distance $0\leq r\leq2$ and an angle $0\leq\alpha\leq\pi$. As shown in
Figure 5, there are two distinct cases to consider.
%TCIMACRO{\FRAME{ftbpFU}{6.1307in}{2.8522in}{0pt}{\Qcb{Two cases for the
%disk.}}{}{las3.eps}{\special{ language "Scientific Word";  type "GRAPHIC";
%width 6.1307in;  height 2.8522in;  depth 0pt;  original-width 6.0701in;
%original-height 2.8098in;  cropleft "0";  croptop "1";  cropright "1";
%cropbottom "0";  filename 'las3.eps';file-properties "XNPEU";}} }%
%BeginExpansion
\begin{figure}[ptb]%
\centering
\includegraphics[
height=2.8522in,
width=6.1307in
]%
{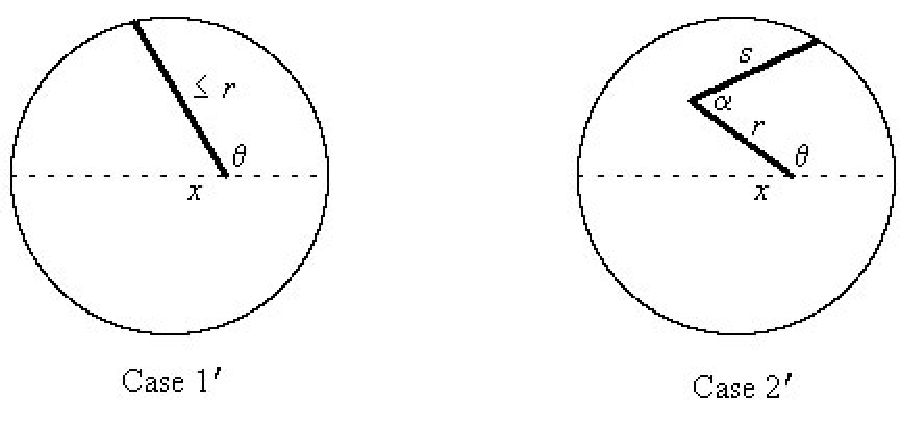}%
\caption{Two cases for the disk.}%
\end{figure}
%EndExpansion
For convenience, define
\[%
\begin{array}
[c]{ccc}%
\varphi(x,r)=\arccos\left(  \dfrac{1-x^{2}-r^{2}}{2xr}\right)  , &  &
\psi(x,r)=\arccos\left(  -\dfrac{r}{x}\right)  .
\end{array}
\]

\textbf{Case 1}$^{\prime}$:\ the swimmer reaches the shore before traveling a
distance $r$, which happens when $x\geq|r-1|$ and $-\varphi\leq\theta
\leq\varphi$, or when $x<r-1$. The path length is $q=-x\cos(\theta
)+\sqrt{1-x^{2}\sin^{2}(\theta)}$.

\textbf{Case 2}$^{\prime}$: the swimmer travels a distance $r$, does not reach
the shore, pivots an angle $\alpha$, and then travels until the shore is
reached. The path length is $r+s$, where
\[
s=y\cos(\alpha+\omega)+\sqrt{1-y^{2}\sin^{2}(\alpha+\omega)},
\]
\[
y=\sqrt{x^{2}+r^{2}+2\,x\,r\cos(\theta)},
\]
\[
\omega=\left\{
\begin{array}
[c]{lll}%
\arcsin\left(  \dfrac{x\sin(\theta)}y\right)  &  & \text{if }\left(  x\geq
r\text{ and }-\psi\leq\theta\leq\psi\right)  \text{ or }x<r,\\
\pi-\arcsin\left(  \dfrac{x\sin(\theta)}y\right)  &  & \text{if }x\geq r\text{
and }\left(  \theta\leq-\psi\text{ or }\psi\leq\theta\right)  .
\end{array}
\right.
\]
The proof of Case 1$^{\prime}$ makes use of the Law of Cosines and the fact
that $\varphi$ exists if and only if $x\geq|r-1|$. The Law of Sines appears in
the proof of Case 2$^{\prime}$ when defining the angle $\omega$ of an
auxiliary triangle with side lengths $r$, $x$ and $y$; the Law of Cosines is
then applied to the triangle with side lengths $1$, $s$ and $y$ to obtain $s$.

The expected path length is therefore $I+J$, where
\[%
\begin{array}
[c]{ccc}%
I=\dfrac{1}{\pi}%
%TCIMACRO{\dint \limits_{|r-1|}^{1}}%
%BeginExpansion
{\displaystyle\int\limits_{|r-1|}^{1}}
%EndExpansion
f\,(x,r,\alpha)\,x\,dx, &  & J=\left\{
\begin{array}
[c]{ccc}%
\dfrac{1}{\pi}%
%TCIMACRO{\dint \limits_{0}^{r-1}}%
%BeginExpansion
{\displaystyle\int\limits_{0}^{r-1}}
%EndExpansion
g\,(x)\,x\,dx &  & \text{if }r\geq1,\\
\dfrac{1}{\pi}%
%TCIMACRO{\dint \limits_{0}^{-r+1}}%
%BeginExpansion
{\displaystyle\int\limits_{0}^{-r+1}}
%EndExpansion
h\,(x,r,\alpha)\,x\,dx &  & \text{if }r<1
\end{array}
\right.
\end{array}
\]
and
\[%
\begin{array}
[c]{ccccc}%
f=%
%TCIMACRO{\dint \limits_{-\varphi}^{\varphi}}%
%BeginExpansion
{\displaystyle\int\limits_{-\varphi}^{\varphi}}
%EndExpansion
q\,d\theta+%
%TCIMACRO{\dint \limits_{\varphi}^{\pi}}%
%BeginExpansion
{\displaystyle\int\limits_{\varphi}^{\pi}}
%EndExpansion
(r+s)\,d\theta+%
%TCIMACRO{\dint \limits_{-\pi}^{-\varphi}}%
%BeginExpansion
{\displaystyle\int\limits_{-\pi}^{-\varphi}}
%EndExpansion
(r+s)\,d\theta, &  & g=%
%TCIMACRO{\dint \limits_{-\pi}^{\pi}}%
%BeginExpansion
{\displaystyle\int\limits_{-\pi}^{\pi}}
%EndExpansion
q\,d\theta, &  & h=%
%TCIMACRO{\dint \limits_{-\pi}^{\pi}}%
%BeginExpansion
{\displaystyle\int\limits_{-\pi}^{\pi}}
%EndExpansion
(r+s)\,d\theta.
\end{array}
\]
Closed-form expressions are not possible here. Numerical computations confirm
that $8/(3\pi)$ is the value of $I+J$ when $\alpha=\pi$, for any $r\geq0$, and
that otherwise $I+J$ is strictly larger than $8/(3\pi)$. See Figure 6 for
several sample realizations (with fixed $r$ but variable $\alpha$).%
%TCIMACRO{\FRAME{ftbpFU}{4.0136in}{4.0318in}{0pt}{\Qcb{Sample realizations with
%$r=0.5$, $x=0.1,0.3,0.5$, $\alpha=110^{\circ},160^{\circ},20^{\circ}$ and
%$\theta=200^{\circ},80^{\circ},45^{\circ}$(all from Case 2$^{\prime}$).}}%
%{}{pth1.eps}{\special{ language "Scientific Word";  type "GRAPHIC";
%width 4.0136in;  height 4.0318in;  depth 0pt;  original-width 3.9652in;
%original-height 3.9825in;  cropleft "0";  croptop "1";  cropright "1";
%cropbottom "0";  filename 'pth1.eps';file-properties "XNPEU";}} }%
%BeginExpansion
\begin{figure}[ptb]%
\centering
\includegraphics[
height=4.0318in,
width=4.0136in
]%
{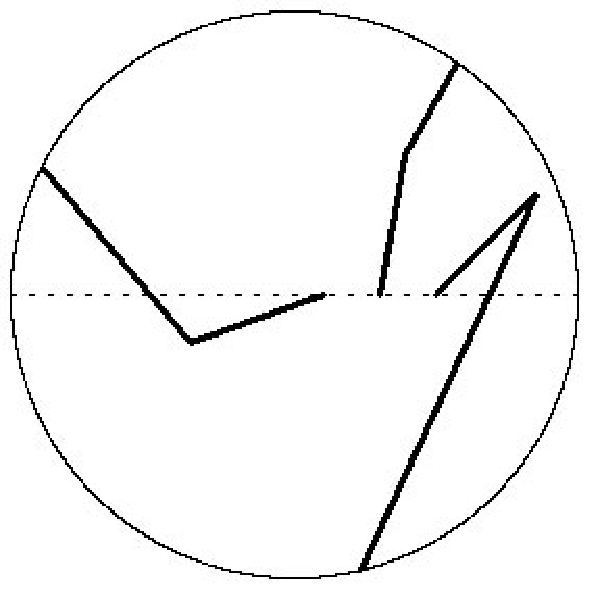}%
\caption{Sample realizations with $r=0.5$, $x=0.1,0.3,0.5$, $\alpha
=110^{\circ},160^{\circ},20^{\circ}$ and $\theta=200^{\circ},80^{\circ
},45^{\circ}$(all from Case 2$^{\prime}$).}%
\end{figure}
%EndExpansion

We mention a corresponding min-median result: Monte Carlo simulation suggests
that a solution is $r=2$ ($\alpha$ can be arbitrary). The path length is the
median of the random variable $-x\cos(\theta)+\sqrt{1-x^{2}\sin(\theta)^{2}}$,
which is estimated to be $0.94...$.

These kinds of ``lost at sea'' problems deserve to be more thoroughly studied!
The only other published reference we know regarding the min-mean problem is
\cite{Gl}.

\subsection{Disk: Gevirtz's Technique}

Let $\gamma:[0,\infty)\rightarrow\mathbb{C}$ be a differentiable curve in the
complex plane satisfying $\gamma(0)=0$. Parametrize $\gamma$ by arclength $s$
and let $\gamma^{\prime}(s)=\exp(i\,\phi(s))$, where we assume that
$\phi(s)=\arg(\gamma^{\prime}(s))$ is real, $\phi(0)=0$ and $|\phi(s)|\leq
\phi_{0}<\frac\pi2$. The function $s\mapsto\lambda(s)=|\gamma(s)|$ is
one-to-one and hence possesses an inverse $\lambda\mapsto s(\lambda)$. Let $D$
denote the unit disk in $\mathbb{C}$ and define $C(z)$ to be the circle of
radius $1$ centered at $z$; as an example, $C(0)$ is the boundary of $D$.

We wish to minimize the average
\[
A(\gamma)=\frac1\pi%
%TCIMACRO{\dint \limits_{D}}%
%BeginExpansion
{\displaystyle\int\limits_{D}}
%EndExpansion
\sigma(z,\gamma)\,dx\,dy
\]
where $\sigma(z,\gamma)$ denotes the first (and only) value of $s$ for which
$|\gamma(s)-z|=1$. (By symmetry of $D$, we have substituted $-z$ for $z$.)
Now,
\[%
%TCIMACRO{\dint \limits_{D}}%
%BeginExpansion
{\displaystyle\int\limits_{D}}
%EndExpansion
\sigma(z,\gamma)\,dx\,dy=%
%TCIMACRO{\dint \limits_{0}^{\infty}}%
%BeginExpansion
{\displaystyle\int\limits_{0}^{\infty}}
%EndExpansion
s\,d\mu(s)
\]
where $\mu(s)$ is the area of the portion of $D$ not enclosed by
$C(\gamma(s))$. This identity follows via a Riemann sum argument.

Consider the family of circular arcs $F=\left\{  D\cap C(\gamma(s)):0<s<\infty
\right\}  $. If $\phi_{0}$ is small enough, then no two distinct arcs in $F$
can intersect. As a consequence of the above identity, if the curvature
$\phi^{\prime}(s)$ is suitably small, then
\[
A(\gamma)=\frac2\pi%
%TCIMACRO{\dint \limits_{0}^{s^{*}}}%
%BeginExpansion
{\displaystyle\int\limits_{0}^{s^{*}}}
%EndExpansion
s\sqrt{1-\left(  \frac{|\gamma(s)|}2\right)  ^{2}}\cos\left(  \phi
(s)-\arg(\gamma(s))\right)  ds
\]
where $s^{*}$ satisfies $|\gamma(s^{*})|=2$. Note that $s^{*}$ depends on
$\gamma$, which complicates any variational approach to this problem. Because
$\lambda^{\prime}(s)=\cos\left(  \phi(s)-\arg(\gamma(s))\right)  $ and
$s(\lambda)\geq\lambda$, we have
\[
A(\gamma)=\frac2\pi%
%TCIMACRO{\dint \limits_{0}^{2}}%
%BeginExpansion
{\displaystyle\int\limits_{0}^{2}}
%EndExpansion
s(\lambda)\sqrt{1-\left(  \frac\lambda2\right)  ^{2}}d\lambda\geq\frac2\pi%
%TCIMACRO{\dint \limits_{0}^{2}}%
%BeginExpansion
{\displaystyle\int\limits_{0}^{2}}
%EndExpansion
\lambda\sqrt{1-\frac{\lambda^{2}}4}d\lambda=\frac8{3\pi},
\]
as was to be shown. It is unclear what to do in the general situation for
which distinct arcs in $F$ are not necessarily disjoint. We also wonder if a
similar technique exists for the infinite strip sea.

\subsection{Acknowledgements}

We are grateful to Julian Gevirtz and Patrice Le Conte for their help.

\end{document}